\documentclass[12pt]{article}
\usepackage{fullpage,graphicx,psfrag,url}
\usepackage[small,bf]{caption}
\usepackage{amsmath,amssymb,amsthm,enumitem}
\usepackage{multirow}
\usepackage{tikz}
\usetikzlibrary{calc}
\newcommand{\BEAS}{\begin{eqnarray*}}
\newcommand{\EEAS}{\end{eqnarray*}}
\newcommand{\BEQ}{\begin{equation}}
\newcommand{\EEQ}{\end{equation}}
\newcommand{\BIT}{\begin{itemize}}
\newcommand{\EIT}{\end{itemize}}

\newcommand{\eg}{{\it e.g.}}
\newcommand{\ie}{{\it i.e.}}

\newcommand{\reals}{{\mbox{\bf R}}}

\newcommand{\symm}{{\mbox{\bf S}}}  

\newcommand{\Tr}{\mathop{\bf Tr}}
\newcommand{\diag}{\mathop{\bf diag}}

\newcommand{\chol}{\mathop{\bf chol}}
\newcommand{\chtwo}{\mathop{\bf chol2}}

\newcounter{oursection}

\newcounter{algorithmctr}[section]
\renewcommand{\thealgorithmctr}{\thesection.\arabic{algorithmctr}}

\title{\Large{Sparsity-Preserving Difference of Positive Semidefinite Matrix Representation of Indefinite Matrices}}
\author{Jaehyun Park}

\begin{document}
\maketitle

\begin{abstract}
We consider the problem of writing an arbitrary symmetric matrix as the
difference of two positive semidefinite matrices. We start with simple ideas
such as eigenvalue decomposition. Then, we develop a simple adaptation of the
Cholesky that returns a difference-of-Cholesky representation of indefinite
matrices. Heuristics that promote sparsity can be applied directly to this
modification.
\end{abstract}

\section{Introduction}

Let $A \in \symm^n$ be a real-valued symmetric matrix.
It is always possible to write $A$
as the difference $P_+ - P_-$ of two positive semidefinite matrices $P_+$ and
$P_-$. Take, for example, $P_+ = A + tI$ and $P_- = tI$ for large enough
$t > 0$. With this representation, one can rewrite an indefinite quadratic
form $x^T A x$ as a difference-of-convex expression $x^T P_+ x - x^T P_- x$,
which can be used in various algorithms for convex-concave
programming~\cite{lipp2014variations,shen2016disciplined}.
In this note, we explore a
number of different such representations and their properties.

\subsection{Desired properties}

In addition to the running time of the algorithm, there are several other properties that are important in typical applications.

\paragraph{Time complexity.}
Clearly, we want the algorithm for computing the representation to be fast.
Asymptotic performance is important, but in practice, even with the same
asymptotic complexity, one method may outperform another depending on how
sparsity is exploited.

\paragraph{Memory usage.}
Ideally, representing $P_+$ and $P_-$ should not cost much more than
representing $P$ itself. In the worst case, both $P_+$ and $P_-$ are
dense, and $2n^2$ floating point numbers must be stored. It is possible to save
more space by exploting the sparsity pattern of $A$, or
representing the matrices $P_+$ and $P_-$ implicitly, in a memory-efficient form.

\paragraph{Numerical stability.}
Due to the round-off errors of floating point numbers, the resulting
representation may be numerically inaccurate. This can be particularly
problematic when $P_+$ and $P_-$ are stored implicitly
(\eg, via Cholesky factorization),
rather than explicitly. Singular or badly conditioned matrices can often
introduce big round-off errors, and we want the algorithm to be robust in such
settings.

\paragraph{Additional curvature.}
In order for convex-concave procedure to perform better, we want $P_+$ and
$P_-$ to have as little ``distortion'' as possible, \ie, we want the
additional curvature measured by the nuclear norm
\[
\|P_+\|_* + \|P_-\|_* - \|A\|_*
\]
to be small. Since $P_+$ and $P_-$ are required to be positive semidefinite, this quantity is equal to
\[
\Tr P_+ + \Tr P_- - \|A\|_*.
\]

\section{Simple representations}

Any indefinite matrix can be made positive semidefinite by adding a large
enough multiple of the identity.
Let $\lambda_{\min} < 0$ be the smallest eigenvalue of $A$. Then, for any
$t \ge |\lambda_{\min}|$,
\[
P_+ = A + t I, \qquad P_- = t I
\]
is a pair of positive semidefinite matrices whose difference is $A$.
If the magnitude of the maximum eigenvalue $\lambda_{\max} > 0$ is smaller than that of $\lambda_{\min}$, then an alternative representation is also possible:
\[
P_+ = tI, \qquad P_- = tI - A,
\]
where $t \ge \lambda_{\max}$. This representation is relatively easy to compute as it only requires a lower bound
on $|\lambda_{\min}|$ or $|\lambda_{\max}|$.
It also has a property that the sparsity of $A$ is
preserved as much as possible, in that no new off-diagonal nonzero entries are
introduced in $P_+$ or $P_-$. Its disadvantage is the additional curvature it introduces:
\[
\|P_+\|_* + \|P_-\|_* - \|A\|_* = 2t.
\]

Another simple representation is based on the full eigenvalue decomposition of
$A$. This representation preserves the norm of $A$ and thus introduces no
additional curvature. Let $A = Q\Lambda Q^T$ be the eigenvalue decomposition
of $A$, where $\Lambda = \diag(\lambda_1, \ldots, \lambda_n)$ is the
eigenvalue matrix with
\[
\lambda_1 \ge \cdots \ge \lambda_k \ge 0 > \lambda_{k+1} \ge \cdots \ge \lambda_n.
\]
Then, $\Lambda$ can be
written as $\Lambda = \Lambda_+ - \Lambda_-$, where
\[
\Lambda_+ = \diag(\lambda_1, \ldots, \lambda_k, 0, \ldots, 0), \qquad
\Lambda_- = \diag(0, \ldots, 0, -\lambda_{k+1}, \ldots, -\lambda_n).
\]
Setting
\[
P_+ = Q\Lambda_+ Q^T, \qquad P_+ = Q\Lambda_- Q^T
\]
then gives a difference of positive semidefinite matrix representation of $A$.
In general, the cost of explicitly computing this representation is very high,
as it not only requires the full eigenvalue decomposition, but also destroys
the sparsity of $A$ even when it is very sparse.

\section{Cholesky-like representations}

When $A$ is positive semidefinite, there exists a unique
lower triangular matrix $L \in \reals^{n \times n}$ that satisfies
$A = LL^T$. This representation is known as the \emph{Cholesky factorization}
of $A$. In the simplest form, the Cholesky algorithm can be described as a
simple recursive algorithm. Formally, let
$\chol: \symm^n \rightarrow \reals^{n \times n}$ that computes,
given $A \in \symm^n$, a lower triangular matrix $L = \chol(A)$ that
satisfies $A = LL^T$. If $A$ is $1$-by-$1$, then $L$ is
simply given by $\sqrt{A_{11}}$. If $A$ has two or more rows, let
\[
A = \left[ \begin{array}{cc} a & v^T \\ v & M \end{array} \right]
\]
with $a \ge 0$, $v \in \reals^{n-1}$, $M \in \reals^{(n-1)\times(n-1)}$.
If $a = 0$ is zero, then $v$ must also equal zero, for otherwise $A$ is 
indefinite. In this case, the Cholesky factorization is given by
\[
\chol \left( \left[ \begin{array}{cc} 0 & 0 \\ 0 & M \end{array} \right] \right)
= \left[ \begin{array}{cc} 0 & \\ 0 & \chol(M) \end{array} \right].
\]
If $a > 0$, the recursion is:
\[
\chol \left( \left[ \begin{array}{cc} a & v^T \\ v & M \end{array} \right] \right)
= \left[ \begin{array}{cc} \sqrt{a} & \\ v/\sqrt{a} & \chol(M - vv^T/a) \end{array} \right].
\]
The cost of computing a dense Cholesky factorization is $(1/3)n^3$ flops. In
case $A$ is sparse, various pivoting heuristics can be used to exploit the
sparsity structure and speed up the computation.

Cholesky factorization does not exist when $A$ has both positive and negative
eigenvalues. However, the LDL decomposition, which is a close variant of the
Cholesky factorization, exists for all symmetric matrices. It also has an
additional computational advantage since there is no need to take square
roots. The idea is to write $A$ as $LDL^T$, where $L \in \reals^{n \times n}$
is lower triangular with ones on the main diagonal, and
$D \in \reals^{n \times n}$ is block diagonal consisting of $1$-by-$1$ or
$2$-by-$2$ blocks. In general, some additional computation is required to
transform the LDL decomposition into the difference of two positive
semidefinite matrices due to the $2$-by-$2$ blocks in $D$. When $D$
is diagonal (\ie, no $2$-by-$2$ blocks), then one can easily separate out the
negative entries of $D$ and the corresponding columns in $L$ to write $A$ as
\[
A = L_1 D_1 L_1^T - L_2 D_2 L_2^T,
\]
where $D_1$ and $D_2$ are diagonal matrices consisting of nonnegative entries
only.

\section{Difference-of-Cholesky representation}

In this section, we describe a simple modification
$\chtwo: \symm^n \rightarrow \reals^{n \times n} \times \reals^{n \times n}$
of the Cholesky algorithm that computes, given $A \in \symm^n$, a pair of
lower triangular matrices $L_1$ and $L_2$ such that
\[
A = L_1 L_1^T - L_2 L_2^T,
\]
where $(L_1, L_2) = \chtwo(A)$. Let $a \ge 0$, $v \in \reals^{n-1}$, and
$M \in \symm^{n-1}$. Below is a simple recursion for computing $\chtwo$
when $a > 0$ is large enough:
\BEAS
\chtwo \left( \left[ \begin{array}{cc} a & v^T \\ v & M \end{array} \right] \right)
&=& \left( \left[ \begin{array}{cc} \sqrt{a} & \\ v/\sqrt{a} & S \end{array} \right],
\left[ \begin{array}{cc} 0 & \\ 0 & T \end{array} \right] \right) \quad (S, T) = \chtwo(M-(1/a) vv^T), \\
\chtwo \left( \left[ \begin{array}{cc} -a & v^T \\ v & M \end{array} \right] \right)
&=& \left( \left[ \begin{array}{cc} 0 & \\ 0 & S \end{array} \right],
\left[ \begin{array}{cc} \sqrt{a} & \\ -v/\sqrt{a} & T \end{array} \right] \right) \quad (S, T) = \chtwo(M+(1/a) vv^T).
\EEAS

In this form, the algorithm is no different from computing the LDL
decomposition and separating out positive and negative entries in $D$
(assuming the LDL decomposition exists with a diagonal $D$). In case
$a$ is too small so that dividing by $\sqrt{a}$ would give big round-off
errors, the following recursion can be used, for any value of $\delta > 0$.
Notice that with this recursion, the algorithm can proceed even when $a=0$.
\[
\chtwo \left( \left[ \begin{array}{cc} a & v^T \\ v & M \end{array} \right] \right)
= \left( \left[ \begin{array}{cc} \sqrt{\delta+a} & \\ v_1 & S \end{array} \right],
\left[ \begin{array}{cc} \sqrt{\delta} & \\ v_2 & T \end{array} \right] \right) \quad (S, T) = \chtwo(M-v_1 v_1^T + v_2 v_2^T).
\]
Here, $v_1$ and $v_2$ are any vectors satisfying
$\sqrt{\delta + a} \, v_1 - \sqrt{\delta} \, v_2 = v$.
This introduces an additional degree of freedom, as we can freely choose $v_1$
or $v_2$. For example, we can let $v_1 = 0$ or $v_2 = 0$.
This way, it is possible to trade off the number
of nonzero elements in $L_1$ against $L_2$, which is a property that is not
readily available in other representations discussed above.
When we choose $v_1 = 0$, the recursion becomes
\[
\chtwo \left( \left[ \begin{array}{cc} a & v^T \\ v & M \end{array} \right] \right)
= \left( \left[ \begin{array}{cc} \sqrt{\delta+a} & \\ 0 & S \end{array} \right],
\left[ \begin{array}{cc} \sqrt{\delta} & \\ -v/\sqrt{\delta} & T \end{array} \right] \right) \quad (S, T) = \chtwo(M+(1/\delta) vv^T).
\]
For completeness, we show the recursion in the case of $A_{11} \le 0$.
\[
\chtwo \left( \left[ \begin{array}{cc} -a & v^T \\ v & M \end{array} \right] \right)
= \left( \left[ \begin{array}{cc} \sqrt{\delta} & \\ v_1 & S \end{array} \right],
\left[ \begin{array}{cc} \sqrt{\delta+a} & \\ v_2 & T \end{array} \right] \right) \quad (S, T) = \chtwo(M-v_1 v_1^T + v_2 v_2^T),
\]
where $v_1$ and $v_2$ are any vectors satisfying
$\sqrt{\delta} \, v_1 - \sqrt{\delta+a} \, v_2 = v$.

It is simple to modify this algorithm to return a pair of LDL factorizations
to avoid computing square roots. With the additional parameter $\delta > 0$,
both factorizations become positive definite and thus there will not be
any $2$-by-$2$ diagonal block in the factorizations.

Since the method is very close to the original Cholesky algorithm, and thus
most extensions and techniques applicable to Cholesky algorithm or LDL
decomposition naturally apply to the difference-of-Cholesky method. One such
example is pivoting heuristics for LDL decomposition for avoiding round-off
errors or ill-conditioned matrices, and at the same time, preserving the
sparsity of the intermediate matrices as much as possible.

\nocite{*}\bibliography{refs}
\end{document}